\documentclass[12pt]{article}

\usepackage{amsmath,amssymb,amsthm,amsfonts}
\usepackage[a4paper,scale=.75]{geometry}
\usepackage{graphicx}
\usepackage{xcolor}
\usepackage{tikz}
\usepackage{bm}
\usepackage{caption}

\newcommand{\N}{\mathbb{N}}
\newcommand{\R}{\mathbb{R}}
\newcommand{\Z}{\mathbb{Z}}
\newcommand{\ga}{\alpha}

\newcommand{\gl}{\lambda}
\newcommand{\gth}{\theta}

\newtheorem{theorem}{Theorem}
\newtheorem{lemma}[theorem]{Lemma}

\newtheorem{proposition}[theorem]{Proposition}
\newtheorem{defn}[theorem]{Definition}

\begin{document}

\title{On Waldron Interpolation on a Simplex in $\R^d$}

\author{Len Bos\and Sione Ma'u\and Shayne Waldron}

\date{}

\maketitle

\begin{abstract}
We introduce explicit families of good interpolation points for interpolation on a triangle in $\R^2$ that may be used for either polynomial interpolation or a certain rational interpolation for which we give explicit formulas. 
\end{abstract}


\section{Introduction}
Although what we claim that our proposed points are very good for interpolation in two dimensions, they may be formulated in all dimensions, and for the sake of generality we do so.  Indeed,  our construction is modelled on (and indeed includes) the so called Simplex Points for a simplex in $\R^d$ which we begin by describing.  

Suppose that ${\bf V}_i\in \R^d,$ $1\le i\le d+1$ are the vertices of a non-degenerate simplex. $S_d.$  For the mulit-index ${\bm \alpha}\in \N_0^{d+1}$ we let
\[{\bm \alpha}=(\alpha_1,\alpha_2,\cdots,\alpha_{d+1}),\quad |{\bm \alpha}|:=\sum_{i=1}^{d+1}\alpha_i.\]
\begin{defn}
For a degree $n\ge 0,$ the Simplex Points of degree $n$ for the simplex $S_d$ are
\[ X_n:=\bigl\{{\bf x}_{\bm \alpha}:=\sum_{i=1}^{d+1}\frac{\alpha_i}{n}{\bf V}_i\,:\, |{\bm \alpha}|=n\bigr\}.\]
\end{defn}
One easily confirms that
\[\#(X_n)={n+d\choose d}={\rm dim}(\Pi_n[\R^d])=:N_n,\]
the polynomials of degree at most $n$ in $d$ real variables. It can also be readily verified that the set of Simplex Points is unisolvent for interpolation by polynomials of degree $n,$ i.e., for every set of
$N_n$ values, $y_{\bm \alpha}\in\R,$ $|{\bm \alpha}|=n,$ there exists a unique polynomial $p\in \Pi_n[\R^d]$ such that
\[p({\bf x}_{\bm\alpha})=y_{\bm \alpha},\quad |{\bm \alpha}|=n.\]
In particular there are explicit formulas (cf.  \cite{Bos:83}) for the cardinal or Lagrange polynomials in barycentric coordinates ${\bm \lambda}\in \R^{d+1},$ $\sum_{i=1}^{d+1}\lambda_i=1:$
\[\ell_{\alpha}({\bm \lambda})= C_{\bm \alpha}\prod_{i=1}^{d+1}\prod_{j=0}^{\alpha_i-1}(\lambda_i-j/n)\]
where $C_{\bm \alpha}$ is a normalization constant given by
\begin{align*}
C_{\bm \alpha}^{-1}&:=\prod_{i=1}^{d+1}\prod_{j=0}^{\alpha_i-1}\Bigl(\frac{\alpha_i}{n}-\frac{j}{n}\Bigr)\cr
&=n^{-n}\prod_{i=1}^{d+1}\alpha_i!.
\end{align*}
The polynomial interpolant may then be written (in barycentric coordinates) in Lagrange form as
\[p({\bm \lambda})=\sum_{|{\bm \alpha}|=n} y_{\bm \alpha}\ell_{\bm \alpha}({\bm \lambda}).\]

The idea behind the interpolation we introduce in this paper (which we refer to as Waldron interpolation),  is to replace the barycentric coordinates ${\bm \alpha}/n$ by a suitable {\it weight} function, $\alpha_i/n \to w(\alpha_i/n).$ However, in general, $\sum_{i=1}^{d+1}w(\alpha_i/n)\neq 1$  and so these are not, strictly speaking, barycentric coordinates.  However, if their  sum is {\it less than} 1, we may create true barycntric coordinates by adding the defect
$(1-\sum_{i=1}^{d+1}w(\alpha_i/n)/(d+1$ to each coordinate, i.e. , setting 
\begin{equation}\label{lambdai}
\lambda_i:= w(\alpha_i/n)+ \frac{1-\sum_{i=1}^{d+1}w(\alpha_i/n)}{d+1},\quad 1\le i\le d+1,
\end{equation}
and consider points expressed as
\[\sum_{i=1}^{d+1} \lambda_i {\bf V}_i,\]
with $\lambda_i$ defined as above.

Notice that if the simplex is centred at the origin, i.e.,
\[\sum_{i=1}^{d+1}{\bf V}_i={\bf 0}\]
then
\[\sum_{i=1}^{d+1}\lambda_i {\bf V}_i=\sum_{i=1}^{d+1}w(\alpha_i/n){\bf V}_i\]
and so the formulas for such points simplify. We will exploit this fact in the sequel.

An immediate example of a possible weight $w$ is $w(x)=x,$ which reduces to the simplex point case.

We will consider  a more general class of weight functions.

\begin{defn} Suppose that $w \,:\,[0,1]\to[0,1].$  We will say that $w$ is an {\bf allowable} weight function if
\begin{itemize}
\item it is increasing on $[0,1],$
\item $w(0)=0$ and $w(1)=1,$ and
\item $\displaystyle{\sum_{i=1}^{d+1}w(\theta_i)\le1}$ for 
$\theta_i\ge0, $ $1\le i\le d+1$ such that $\displaystyle{\sum_{i=1}^{d+1}\theta_i=1}.$
\end{itemize}
We remark that if $w_0$ and $w_1$ are two allowable weight functions, then the convex combinations
\[w_t:=tw_1+(1-t)w_0,\quad 0\le t\le 1,\]
are also allowable.
\end{defn}

\begin{proposition}\label{allowable}
Let $F$ be a nonnegative, non-decreasing function on $[0,1/2]$, which 
is normalised so that
\[\int_0^{1/2} F(t) = \frac{1}{2}.\]
Let $w$ be the  increasing convex continuous function on $[0,1/2]$
with $w(1/2)=1/2$,
defined by 
\[w(x) = \int_0^x F(t)\,dt, \quad 0\le x\le 1/2,\]
and extended to $[0,1]$ by 
\[w(x):=\int_0^x \widetilde{F}(t)\,dt\]
with
\[
\widetilde{F}(t)
=\begin{cases}
F(t), & 0\le t\le 1/2; \cr
F(1-t), & 1/2< t\le1.
\end{cases}
\]
Then $w$ is an allowable function which satisfies
\[w(\gth_1)+\cdots+w(\gth_{d+1}) \le w(\gth_1+\cdots+\gth_{d+1})\le w(1)=1,
\qquad \gth_j\ge0, \quad \sum_{j=1}^{d+1}\gth_j\le1.\]
Furthermore, $w$ is complementary, i.e.,
\[w(x)+w(1-x)=1, \quad 0\le x\le 1, \]
and
\[w(x)\le x, \quad 0\le x\le{1/2}, 
\qquad w(x)\ge x, \quad  {1/2}\le x\le 1. \]
\end{proposition}

\begin{proof}
By construction, $w$ is clearly increasing with $w(0)=0$ and $w(1)=1$.
Thus, to show that $w$ is allowable, it suffices to prove the inequality
\[w(\gth_1)+\cdots+w(\gth_{d+1}) \le w(\gth_1+\cdots+\gth_{d+1}), 
\qquad \gth_j\ge0, \quad \sum_{j=1}^{d+1}\gth_j\le1, \]
by induction on $d\ge1$. We now prove the case $d=2$ (it is trivial for $d=1$),
which is used in the inductive step. 
We observe that
\begin{align*}
w(\gth_1)+w(\gth_2)\le w(\gth_1+\gth_2)
& \iff
\int_0^{\gth_1}\widetilde{F}(t)dt +\int_0^{\gth_2}\widetilde{F}(t)dt
\le \int_0^{\gth_1+\gth_2}\widetilde{F}(t)dt\cr
&\qquad = \int_0^{\gth_1}\widetilde{F}(t)dt + \int_{\gth_1}^{\gth_1+\gth_2}\widetilde{F}(t)dt\cr
& \iff
\int_0^{\gth_2}\tilde{F}(t)dt\le \int_{\gth_1}^{\gth_1+\gth_2}\widetilde{F}(t)dt,
\end{align*}
which holds since $\widetilde{F}$ is increasing. 

Now suppose the result holds for $d-1.$ Then
\begin{align*}
w(\gth_1)+\cdots+w(\gth_{d+1})
&=\{w(\gth_1)+\cdots+w(\gth_d)\} +w(\gth_{d+1}) \cr
&\le w(\gth_1+\cdots+\gth_d)+w(\gth_{d+1}) \cr
&\le w(\{\gth_1+\cdots+\gth_d\}+\gth_{d+1}),
\end{align*}
which completes the induction.

Suppose,  without loss of generality, that $x\le{1\over2}.$ Then
\begin{align*}
w(x)+w(1-x)
&= \int_0^x F(t)\, dt + \left\{ \int_0^{1\over2} F(t)\,dt
+\int_{1\over2}^{1-x} F(1-t)\,dt \right\} \cr
&= \int_0^x F(t)\, dt + {1\over2}
+\int_x^{1\over2} F(s)\,ds \quad (\hbox{letting}\,\,s:=1-x)\cr
&= {1\over2} + \int_0^{1\over2} F(t)\, dt\cr
&={1\over2}+{1\over2}=1.
\end{align*}
The condition that $w$ is convex on $[0,{1\over2}]$ is
$$ w({t\over2})=w((1-t)0+t{1\over2})\le(1-t)w(0)+tw({1\over2})
= {t\over 2},
\qquad 0\le t\le 1. $$
i.e., $w(x)\le x$, $0\le x\le{1\over2}$. 
The inequality for $x\ge{1\over2}$, then follows by the calculation
$$ w(x)=1-w(1-x)\ge1-(1-x)=x, $$
since $0\le1-x\le{1\over2}$.
\end{proof}

\noindent{\bf Example 1}.  If we take $F(t)=1$, then $\widetilde{F}(t)\equiv 1,$
and
\[w(x)= \int_0^x \widetilde{F}(t)dt=x\]
and we recover the Simplex Points.

\medskip

\noindent {\bf Example 2}.  If we take $F(t)=\frac{\pi}{2}\sin(\pi t)$ then
\begin{align*}
\widetilde{F}(t)&=\begin{cases}
 \frac{\pi}{2}\sin(\pi t)& 0\le t\le 1/2; \cr
\frac{\pi}{2}\sin(\pi (1-t))& 1/2< t\le1
\end{cases}\cr
&=\begin{cases}
 \frac{\pi}{2}\sin(\pi t)& 0\le t\le 1/2; \cr
\frac{\pi}{2}\sin(\pi t)& 1/2< t\le1
\end{cases}\cr
&=\frac{\pi}{2}\sin(\pi t), \quad 0\le t\le 1.
\end{align*}
Hence
\begin{align*}
w(x)&=\int_0^x \widetilde{F}(t)dt\cr
&=\frac{\pi}{2}\int_0^x \sin(\pi t)dt\cr
&=\frac{1-\cos(\pi x)}{2} =\sin^2(\pi x/2).
\end{align*}

\medskip

\noindent {\bf Example 3}.  Take $F(t):=4t$ for which it results that
\[w(x)=\left\{\begin{array}{cc}
2x^2,&0\le x\le 1/2\cr
\cr
1-2(1-x)^2,& 1/2\le x\le 1\end{array}\right. .\]
$\Box$

%
%

\medskip

For an allowable weight function we introduce the associated Waldron points:

\begin{defn}\label{waldronpts}
Suppose that $w(x)$ is an allowable weight function.  The associated Waldron points of degree $n$ for a simplex $S_d\subset \R^d$ with vertices ${\bf V}_i,$ $1\le i\le d+1,$ are given by
\[W_n:=\Bigl\{ {\bf x}_{\bm \alpha}=\sum_{j=1}^{d+1} \omega_j {\bf V}_j\,:\,|{\bm \alpha}|=n\Bigr\}\]
where
\[ \omega_j:= w(\alpha_j/n)+\frac{1}{d+1}\Bigl(1-\sum_{i=1}^{d+1}w(\alpha_i/n)\Bigr),\,\,1\le j\le d+1.\]
In the case that $S_d$ is centred at the origin, i.e.,
\[\sum_{j=1}^{d+1}{\bf V}_j={\bm 0}\in\R^d,\]
then we may take
\[\omega_j:= w(\alpha_j/n),\quad 1\le j\le (d+1).\]
We refer to the $\omega_j$ as the {\it baryweights} of the Waldron points.
\end{defn}

We let 
\[T_d:=\{{\bm \theta}\in \R^{d+1}\,:\, \theta_j\ge0,\,\,\sum_{j=1}^{d+1}\theta_j=1\}\]
be the standard unit simplex in the positive orthnant, consisting of all possible barycentric coordinates for points in a general simplex.  We show that in dimension $d=2,$ and a triangle centred at the origin, the baryweights for an allowable weight form a coordinate system on the triangle.

\begin{proposition} \label{Onto}
Consider the simplex $S_d\subset \R^d$ with vertices ${\bf V}_i,$ $1\le i\le d+1,$ centred at ${\bf 0}$, i.e., with $\sum_{j=1}^{d+1}{\bf V}_j=0$. If $w$ is an 
allowable weight, then the map between simplices
\begin{equation}
\label{maponsimplices}
T_d\to S_d: (\gth_j)\mapsto x=\sum_{j=1}^{d+1} w(\gth_j){\bf V}_j,
\end{equation}
is $1$--$1$ with image
\[ S_d^w := \{ {\bf x}=\sum_{j=1}^{d+1}\gl_j{\bf V}_j:{\bm \lambda}\in T_d,\
\sum_{j=1}^{d+1} w^{-1}(\gl_j-\gl_{\rm min})\le1\}\subset S_d; \qquad \gl_{\rm min}=\min_{1\le j\le d+1}\gl_j. \]
Moreover, if $w$ is complementary and $d=2$, then $S_2^w=S_2$, i.e.,
each point ${\bf x}$ in the triangle $S_2$ has unique baryweights $(w(\gth_j))$.
\end{proposition}

\begin{proof} Let ${\bf x}=\sum_{j=1}^{d+1} \gl_j{\bf V}_j$, ${\bm \lambda}\in T_d,$ be an
arbitrary point of $S_d$. 
We consider the condition that there exists a ${\bm \theta}\in T_d$ with 
\[ {\bf x}=\sum_{j=1}^{d+1} w(\gth_j){\bf V}_j. \]
Since the barycentric coordinates are unique, this can be written as
\[ \gl_j = w(\gth_j)+{1\over d+1}\Bigl(1-\sum_{k=1}^{d+1} w(\gth_k)\Bigr), 
\qquad\forall j, \]
which is equivalent to the square linear system
\[ d\, w(\gth_j) -\sum_{k\ne j} w(\gth_k) = (d+1)\gl_j-1, \qquad
1\le j\le (d+1) \]
in the variables $a_j=w(\gth_j)$. For $d\ge1$, the matrix of the (diagonal entries $d$, and off-diagonal entries $-1$) has a $1$-dimensional
kernel spanned by $(1,1,\ldots,1)$, and so if there is a solution $(a_j)$,
then all solutions are given by $(a_j+c)$, $c\in\R$.  One such solution is to  take
 $a_j=\gl_j$,  $1\le j \le d+1.$  Thus all possible
choices for $(\gth_j)$ are given by
\[w(\gth_j)=\gl_j+c \iff \gth_j=w^{-1}(\gl_j+c), \qquad\forall j, \]
for a suitable constant $c$. For the above formula 
to give a point ${\bm \theta}$  in the simplex $T_d$, this $c$ must satisfy
\[ -\min_j\gl_j\le c\le 1-\max_j \gl_j, \qquad
H(c):=\sum_{j=1}^{d+1} w^{-1}(\gl_j+c) = 1.\]
Since $w^{-1}$ is (strictly) increasing and continuous, 
so is $c\mapsto H(c)$, and hence 
there can be at most one such choice for $c$, 
i.e., the map between simplices is $1$--$1$.
By Lemma \ref{wwinvoptimistation}, we have
\[ H(0)=\sum_{j=1}^{d+1}w^{-1}(\gl_j) \ge1, \]
and so, by the intermediate value theorem, 
there is a (unique) choice for $c\le0$
(and hence ${\bf x}$ is mapped onto) if and only if
\[ H(-\min_j\gl_j) = \sum_j w^{-1}(\gl_j-\gl_{\rm min})\le1. \]

Now suppose $w$ is complementary and $d=2$. Since $w$ is complementary,
we have
\[w(x)=1-w(1-x) \implies w^{-1}(1-y)=1-w^{-1}(y). \]
Assume, with out loss of generality, that $\gl_1\le\gl_2\le\gl_3$, then
we calculate
\begin{align*} 
\sum_{j=1}^3w^{-1}(\gl_j-\gl_{\rm min})
& =w^{-1}(0)+w^{-1}(\gl_2-\gl_1)+w^{-1}(\gl_3-\gl_1)\cr
&\le w^{-1}(\gl_2)+w^{-1}(\gl_3) \cr
&\le w^{-1}(1-\gl_3)+w^{-1}(\gl_3)\quad (\hbox{as}\,\, \lambda_2+\lambda_3\le1)\cr
&=( 1-w^{-1}(\gl_3) )+w^{-1}(\gl_3)\cr
&=1,
\end{align*}
which shows $S_2^w=S_2$.
\end{proof}

\begin{lemma}
\label{wwinvoptimistation}
Let $w$ be an allowable weight given by a non-decreasing $F$. Then for ${\bm \theta}\in T_d$, we have
\[ (d+1)w({1\over d+1}) \le \sum_{j=1}^{d+1} w(\gth_j)\le 1, \]
\[ 1 \le \sum_{j=1}^{d+1} w^{-1}(\gth_j)\le (d+1)w^{-1}({1\over d+1}), \]
where there is an equality for ${\bm \theta}$ a vertex or the barycentre of $T_d$.
\end{lemma}
\begin{proof}
We proceed by induction on $d.$ The case $d=1$ is trivial as $w(x)$ is complementary.  Hence assume that the conclusions hold for all dimensions $d'<d.$

We remark for that $y\le 1/2,$ as $F$ is non-decreasing,
\[F(y)\ge\frac{1}{y}\int_0^y F(t)dt.\]
Hence, setting $y=1/x,$
\begin{align*}
\int_0^{1/x}F(t)dt-\frac{1}{x}F(1/x)&\le0\cr
\implies \,\,w(1/x)-\frac{1}{x}w'(1/x)&\le0\cr
\implies \,\, \frac{d}{dx}\Bigl\{x\,w(1/x)\Bigr\}&\le 0
\end{align*}
and so
\[(d+1)w(\frac{1}{d+1})\] 
is {\it decreasing} in the dimension $d.$

Similarly,
\[(d+1)w^{-1}(\frac{1}{d+1})\]
is {\it increasing} in the dimesnion $d.$

It follows that if ${\bm \theta}\in T_d$ is a boundary point, i.e., one of the $\theta_j=0,$ say $\theta_k=0,$ we have that
\begin{align*}
\sum_{j=1}^{d+1}w(\theta_j)&=\sum_{j\neq k}w(\theta_j)\cr
&\ge d\,w(\frac{1}{d})\cr
&\ge (d+1)\,w(\frac{1}{d+1})
\end{align*}
and simlilarly
\begin{align*}
\sum_{j=1}^{d+1}w^{-1}(\theta_j)&=\sum_{j\neq k}w^{-1}(\theta_j)\cr
&\le d\,w^{-1}(\frac{1}{d})\cr
&\le (d+1)\,w^{-1}(\frac{1}{d+1}).
\end{align*}

The upper bound of $1$ for the sum of the $w(\theta_j)$ is a consequence of being an allowable weight (see Proposition
\ref{allowable}).  For the lower bound on the sum of the inverse weight functions,  if ${\bm \theta}\in T_d$ has a zero component,  $\theta_k,$ say, then
\[\sum_{j=1}^{d+1}w^{-1}(\theta_j)=\sum_{j\neq k} w^{-1}(\theta_j)\ge1\]
by the induction assumption. 

Hence we may assume that ${\bm \theta}$ is not a boundary point of $T_d.$

Since ${\bm \theta}$ satisfies $\sum_{j=1}^{d+1}\gth_j=1$,
we can eliminate one variable, say 
\[\gth_1=1-\sum_{j=2}^{d+1}\gth_j,\]
and then seek to optimise
\[ f(\gth_2,\ldots,\gth_{d+1})
:= \sum_{j=1}^{d+1} w(\gth_j)
= w\bigl(1-\sum_{j=2}^{d+1} \gth_j\bigr)+\sum_{j=2}^{d+1} w(\gth_j), \qquad
\gth_2+\cdots+\gth_{d+1}\le1. \]
This has a critical point when
\[ {\partial f\over\partial \gth_k}
= -w'\Bigl(1-\sum_{j=2}^{d+1}\gth_j\Bigr)+ w'(\gth_k)
=0, \quad\forall k, \]
i.e., when $w'(\gth_k)$ is a constant given by
\[ w'(\gth_k) = w'(1-\sum_{j=2}^{d+1} \gth_j), \qquad k=2,\ldots,(d+1). \]
First suppose that $w'=F$ strictly increasing on $[0,{1\over2}]$, so that
$w'$ takes a given value at most twice on $[0,1]$, i.e., at $x$ and $1-x$,
and the condition for a critical point gives, for $k=2,3,\cdots,(d+1),$ either
\[ \gth_k=1-\sum_{j=2}^{d+1}\theta_j ,\]
or
\[ \gth_k=\sum_{j=2}^{d+1} \gth_j.\]
If the second case occurs,  then
\[\theta_k=\theta_k+\sum_{j=2,\,\,j\neq k}^{d+1}\theta_j\,\,\implies \,\, \theta_j=0,\,\,j\neq k\]
and we are at a boundary point.

Hence the first case must always hold.  However, the matrix of this linear system (with $2$ on the diagonal and $1$ on the off-diagonal) is non-singular and easily solved to find that
\[\theta_j=\frac{1}{d+1},\quad 1\le j\le d+1,\]
i.e., ${\bm \theta},$ the barycentre of $T_d$  is the {\it  unique} critical point in the interior of $T_d.$ SInce $f({\bm\theta})=(d+1)w(\frac{1}{d+1})$ and, by the induction hypothesis then,  $f\ge f({\bm\theta})$ on the boundary of $T_d, $ it follows that  the barycentre is a {\it minimum} point and $f({\bf x})\ge (d+1)w(\frac{1}{d+1})$ for all ${\bf x}\in T_d.$

If $w'$ is not strictly increasing, e.g., $w(x)=x$, then take a 
convex combination $w_\ga=(1-\ga)w+\ga w_s$, $0<\ga\le1$ with a $w_s$ which
is strictly increasing, so that $w_\ga$ is strictly increasing, and
therefore satisfies the inequalities. Finally taking the limit of
these inequalities as $\ga\to0^+$ gives the result for $w$.

The proof for $w^{-1}$, which is increasing, is essentially the same,
since the corresponding ``$F$'' on $[0,{1\over2}]$ which is given by
$$ (w^{-1})'(x)={1\over w'(w^{-1})'(x))}, \qquad x\ne0, $$
is decreasing, possibly having an integrable pole at $x=0$
(when $w'(0)=0$). 
\end{proof}

\section{A Geometric Interpretation}

For the interval $[-1,1]$ the extended Chebyshev points
\[x_k=\cos(k\pi/n),\quad 0\le k\le n\]
are near optimal for polynomial interpolation.

In barycentric coordinates they may be expressed as
\begin{align*}
x_k&=\frac{1-\cos(k\pi/n)}{2}\times (-1) + \frac{1+\cos(k\pi/n)}{2}\times (+1)\cr
&= w(k/n)\times (-1)+ w((n-k)/n)\times (+1)
\end{align*}
for 
\[w(x):=\frac{1-\cos(\pi x)}{2}.\]

In several variables the spacing of near optimal points is more complicated.  It is known that a sequence of point sets for which the Lebesgue function is of sub-exponential growth,  the equally weighted discrete probability measures supported on these points must tend weak-* to the so-called equilibrium measure of Pluripotential Theory (cf.  \cite{K}). 

As discussed in \cite{BLW:04}, natural distances capturing the spacing of these point sets may be associated to the equlibrium measure.  In particular, for a simplex in $\R^d,$ the equilibrium measure is, in terms of the barycentric coordinates,
\[ \mu_{S_d}=c \frac{1}{\sqrt{\lambda_1\lambda_2\cdots\lambda_{d+1}}}dA\]
where $c$ is a normalization constant to make it a probability measure.  The so-called associated Baran distance is given as follows. Suppose that ${\bf a},{\bf b}\in S_d$ are two points in the simplex. Then
\begin{equation}\label{BaranDist}
d_B({\bf a},{\bf b}):=\cos^{-1}(\widetilde{\bf a}\cdot \widetilde{\bf b})
\end{equation}
where, for ${\bf x}\in S_d$ with barycentric coordinates $\lambda_1,\cdots,\lambda_{d+1},$
\[\widetilde{\bf x}:=(\sqrt{\lambda_1},\sqrt{\lambda_2},\cdots,\sqrt{\lambda_{d+1}})\in \mathbb{S}_d,\]
the unit sphere in $\R^{d+1}.$

Indeed, the map ${\bf x}\to \widetilde{\bf x}$ maps the simplex onto the positive orthnant of the sphere
\[\{\widetilde{\bf x}\in \R^{d+1}\,:\, \|\widetilde{\bf x}\|_2=1,\,\, \widetilde{x}_j\ge0,\,\,1\le j\le d+1\}\]
with ``vertices'' at 
\[(1,0,\cdots,0),(0,1,0,\cdots,0),\cdots,(0,\cdots,0,1).\]
The geodesics on the sphere,  emanating from a vertex, say from $(1,0,\cdots,0),$ are great circles and may be described in polar-hyperspherical coordinates, as
\[(\cos(\phi_1),\sin(\phi_1){\bf y}), \quad 0\le \phi_1\le \pi/2\]
where ${\bf y}\in\mathbb{S}_{d-1}$ is fixed.

This pulls back to the simplex by
\[\widetilde{\bf x} \to {\bf x}\]
with
\[ \lambda_j=\widetilde{x_j}^2,\quad,1\le j\le d+1\]
i.e.,
\begin{align*}
\lambda_1&=\cos^2(\phi_1)\cr
\lambda_j&=\sin^2(\phi_i)y_j^2,\quad 2\le j\le (d+1)\cr
&=(1-\lambda_1)y_j^2
\end{align*}
which describes a straight line passing through the vertex ${\bf V}_1$ of the simplex.  Moreover, the hyperplanes 
\[\lambda_j=\frac{1-\cos(k\pi/n)}{2},\quad 0\le k\le n\]
are equally spaced in this Baran distance.

It would seem to be ideal if there were points at the intersections of $d+1$ of these hyperplanes, $\lambda_j=w(\alpha_j/n).$ However, in general this is not possible and perhaps a good compromise would be to consider the simplex generated by the $d+1$ intersections of the $d+1 \choose d = d+1$ subsets of $d$ of the $d+1$ hyperplanes, and then take the centroid of these simplices as our interpolation points.  

To see what this entails, note that each combination of $d$ of the $d+1$ hyperplanes amounts to leaving one of the equations $\lambda_j=w(i/n)$ out.  If the $k$th equation is left out then the intersection point (in barycentric coordinates) is
\begin{align*}
\lambda_j&=w(\alpha_j/n),\quad j\neq k,\cr
\lambda_k&=1-\sum_{j\neq k} w(\alpha_j/n).
\end{align*}
Their average is therefore
\begin{align*}
\lambda_j&=\frac{d \,w(\alpha_j/n)+\Bigl(1-\sum_{k\neq j} w(\alpha_k/n)\Bigr)}{d+1}\cr
&=\frac{(d+1)w(\alpha_j/n)+\Bigl(1-\sum_{k=1}^{d+1}w(\alpha_k/n)\Bigr)}{d+1}\cr
&= w(\alpha_j/n)+\frac{1}{d+1}\Bigl(1-\sum_{k=1}^{d+1}w(\alpha_k/n)\Bigr),
\end{align*}
precisely the Waldron points of Definition \ref{waldronpts}.  Figure \ref{fig1} shows a typical example.

\begin{center}
 \includegraphics[width=10cm,height=10cm]{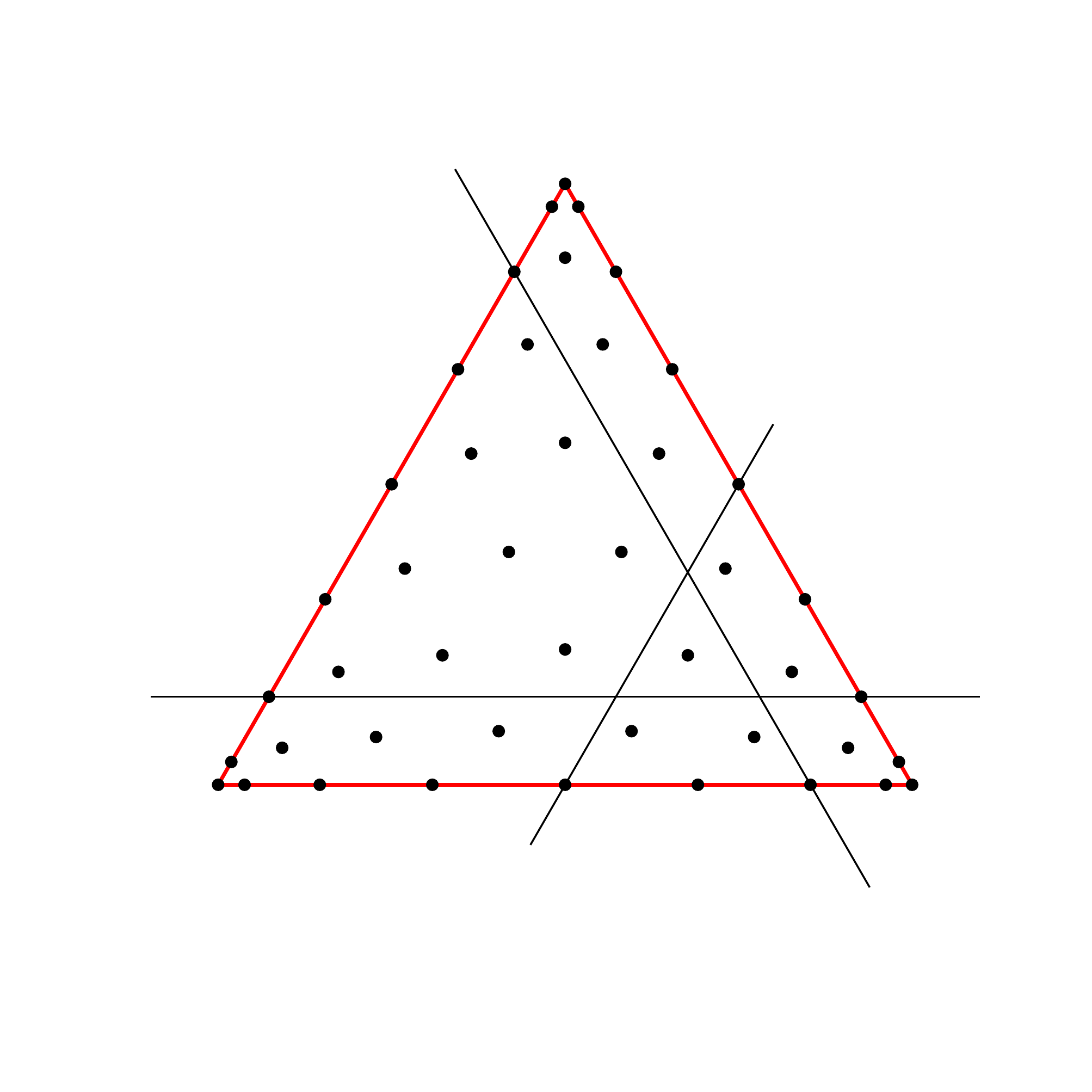}
\captionof{figure}{Triangle Formed by Coordinate Lines}
\label{fig1}
\end{center}

{\bf Remark}.  In the case that $w(x)=x$ and we recover the Simplex Points, the (d+1) coordinate hyperplanes do indeed intersect at a single point.

\begin{center}
 \includegraphics[width=10cm,height=10cm]{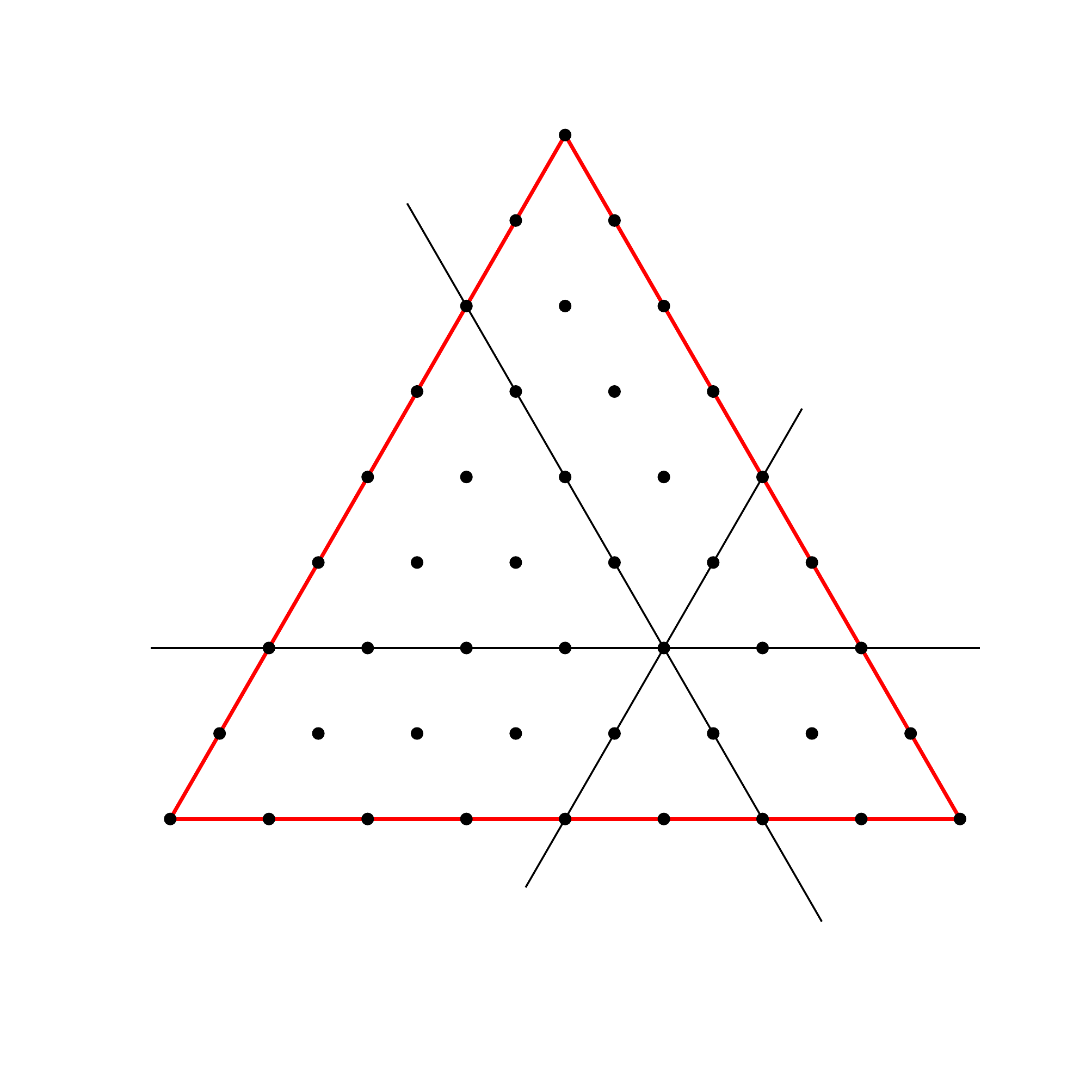}
\captionof{figure}{Coordinate Lines for Simplex Points}
\label{fig2}
\end{center}


\section{A Triangle in $\R^2$ Centred at the Origin}

In this section we consider a triangle $S_2\subset \R^2$ with vertices ${\bf V}_1,{\bf V}_2,{\bf V}_3\in \R^2$ centred at the origin, i.e., with
\[ {\bf V}_1+{\bf V}_2+{\bf V}_3=(0,0).\]
For ${\bm \theta} \in \R^3,$ $\sum_{i=1}^3\theta_i=1,$ $\theta_1,\theta_2,\theta_3\ge0,$ we may use barycentric coordinates, as in (\ref{lambdai}),
\[\lambda_i=w(\theta_i)+\frac{1}{3}\Bigl(1-\sum_{i=1}^3 w(\theta_i)\Bigr),\quad i=1,2,3,\]
to define a point
\begin{align*}
{\bf x}&=\sum_{i=1}^3 \lambda_i {\bf V}_i\cr
&=\sum_{i=1}^3 w(\theta_i) {\bf V}_i\in S_2
\end{align*}
(using the centredness of the vertices).

Now,  if every point ${\bf x}\in S_2$ may be uniquely represented in such a manner,  we have a legitimate coordinate system on the triangle and we refer to the $w(\theta_i)$ as {\it baryweights}.

Note also that in  the case of a simplex centred at the origin, the Waldron points of Definition \ref{waldronpts} are 
\[ W_n:=\Bigl\{ {\bf x}_{\bm \alpha}=\sum_{j=1}^{d+1} w(\alpha_j/n) {\bf V}_j\,:\,|{\bm \alpha}|=n\Bigr\}.\]
In other words, the Waldron points are those with baryweights
$w(\alpha_i/n),$ $|{\bm\alpha}|=n.$

For the Waldron points we may then define cardinal (Lagrange) functions
\begin{equation}\label{newLagrange}\ell_{\bm \alpha}({\bf x})= C_{\bm \alpha}\prod_{i=1}^{3}\prod_{j=0}^{\alpha_i-1}(w(\theta_i)-w(\alpha_j/n))
\end{equation}
for ${\bm \theta}\in \R_+^{3},$ $\sum_{j=1}^{3}\theta_j=1,$ such that
\[{\bf x}=\sum_{j=1}^{3}w(\theta_j)V_j.\]
Here $C_{\bm \alpha}$ is a normalization constant given by
\[
C_{\bm \alpha}^{-1}:=\prod_{i=1}^{3}\prod_{j=0}^{\alpha_i-1}\Bigl(w(\alpha_i/n)-w(\alpha_j/n)\Bigr).
\]

An interpolant may then be written  in Lagrange form as
\begin{equation}\label{newinterp}
q({\bf x})=\sum_{|{\bm \alpha}|=n} y_{\bm \alpha}\ell_{\bm \alpha}({\bf x}).
\end{equation}

We will show some examples of this interpolant below.  But first note that for $w(x)=x$ we recover polynomial interpolation at the Simplex Points.  By Proposition \ref{Onto} our other example of an allowable weight function 
\[w(x)=\frac{1-\cos(\pi x)}{2}\]
also provides a legitimate coordinate system.

\medskip
In particular we may consider the interpolant (\ref{newinterp}) with $w(x)=(1-\cos(\pi x))/2.$ It is {\it not} a polynomial interpolant and indeed does not even reproduce the constant function $f({\bf x})=1.$ Figure \ref{fig5} shows the interpolant of 1 for degree $n=2.$

\begin{center}
 \includegraphics[width=12cm,height=12cm]{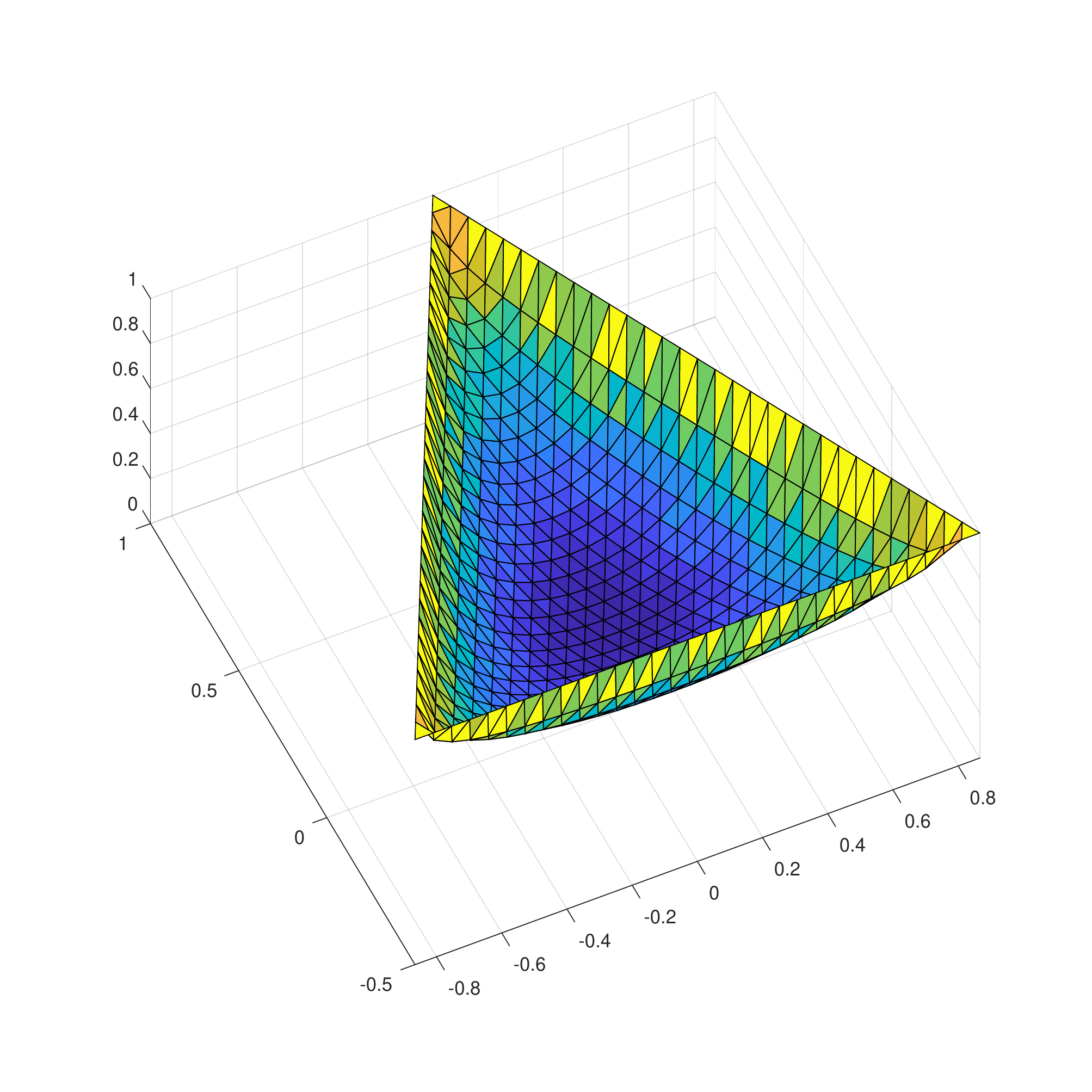}
\captionof{figure}{Interpolant of 1 for $n=2$}
\label{fig5}
\end{center}

However, if we normalize the Lagrange polynomials (\ref{newLagrange}) by dividing by their sum,  
\[\ell_{\bm \alpha}({\bf x}) \to \frac{\ell_{\bm \alpha}({\bf x})}
{\sum_{|{\bm \beta}|=n}\ell_{\bm \beta}({\bf x})},\]
 we obtain a superior interpolant, albeit now a rational function. It retains the distinct advantage of having explict points and explicit formulas for the cardinal functions. Figure \ref{fig6} shows the interpolant of the function $f(x,y)=\sin(\pi(x^2+y^2))$ for degree
 $n=5.$
 
 \begin{center}
 \includegraphics[width=12cm,height=12cm]{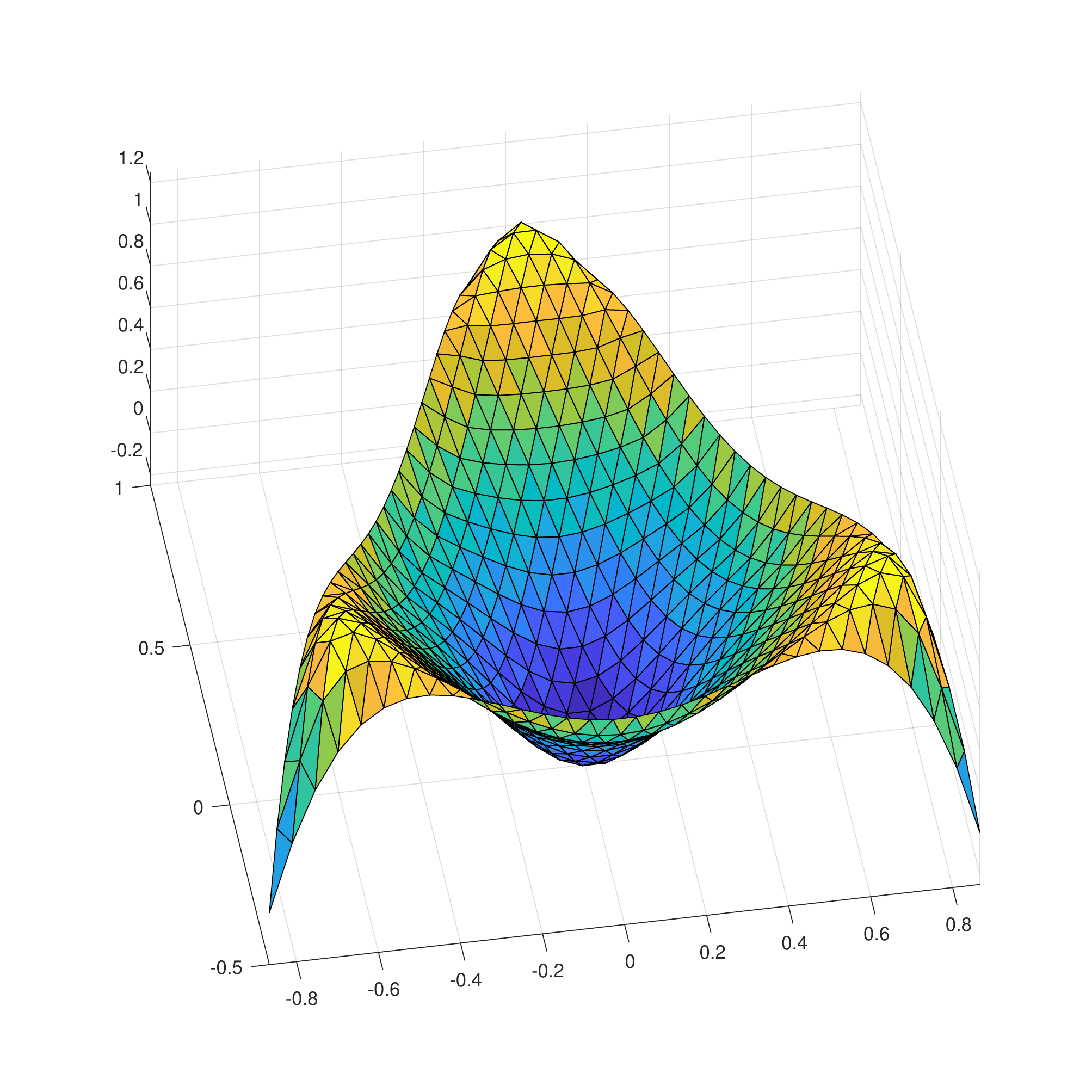}
\captionof{figure}{Interpolant of $f(x,y)$ for $n=5$}
\label{fig6}
\end{center}

\subsection{Spherical Waldron Points}
The map from a triangle to the sphere given by
\[ {\bf x}\to (\sqrt{\lambda_1({\bf x})},\sqrt{\lambda_2({\bf x})},
\sqrt{\lambda_3({\bf x})})\]
maps onto the positive octant $x,y,z\ge0$ and pushes forward the Baran distance on the triangle to spherical distance on the sphere (as is clear from the formula for the Baran distance (\ref{BaranDist})). Hence points, well spaced with respect to the Baran distance on the triangle, will map onto nearly equally spaced points on the sphere. 

Here we show the points for degree $n=20$ that result from taking the square roots of the Waldron weights and then normalizing them to lie on the sphere,i.e, the set of points
\[SW_n:=\Bigl\{ {\bf x}_{\bm \alpha}= \frac{(\sqrt{w(\alpha_1/n)},\sqrt{w(\alpha_2/n)},\sqrt{w(\alpha_3/n)})}{\sqrt{w(\alpha_1/n)+w(\alpha_2/n)+w(\alpha_3/n)}}\,:\,|{\bm \alpha}|=n\Bigr\}\]
where
\[w(x)=\frac{1-\cos(\pi x)}{2}=\sin^2\Bigl(\frac{\pi}{2}x\Bigr).\]

Obviously, these may be combined by rotation and reflection to obtain a set of $4n^2+2$ points on the full sphere.

\begin{center}
 \includegraphics[width=12cm,height=12cm]{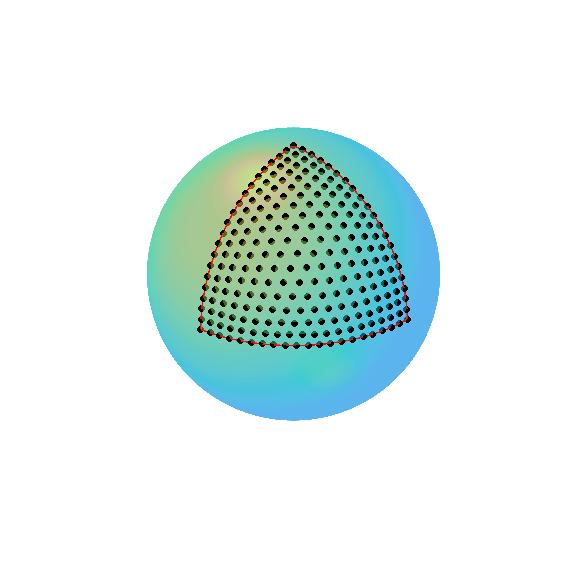}
\captionof{figure}{Spherical Waldron Points for Degree 20}
\label{fig4}
\end{center}

The spacing of these points may be analyzed.  Set
\[v(\theta):=\sqrt{w(\theta)}\,\,\hbox{and}\,\, W({\bm\theta}):=w(\theta_1)+w(\theta_2)+w(\theta_3).\]
Then 
\[X({\bm \theta}):=\Bigl(\frac{v(\theta_1)}{\sqrt{W({\bm \theta})}},
\frac{v(\theta_1)}{\sqrt{W({\bm \theta})}},
\frac{v(\theta_1)}{\sqrt{W({\bm \theta})}}\Bigr)\]
maps $T_2$ onto the sphere in the positive octant, such that the Waldron points are
\[{\bf X}({\bm \alpha}/n),\quad {\bm \alpha}\in \Z_+^3,\,\,|{\bm \alpha}|=n.\]
Neighbouring index points are those whose indices differ by $+1$ in one component and by $-1$ in one of the two others.  Without loss we may consider these two components to be the frist and the second.  Then the (euclidean) distance between two neighbouring Waldron points is
\[\|{\bf X}(\theta_1+1/n,\theta_2-1/n,\theta_3)-{\bf X}(\theta_1,\theta_2,\theta_3)\|_2.\]
Setting $h=1/n$ this may be expressed as
\[h\,\bigl\|\frac{{\bf X}(\theta_1+h,\theta_2-h,\theta_3)-{\bf X}(\theta_1,\theta_2,\theta_3)}{h}\bigr\|_2.\]
This in turn is asymptotically (as $n\to\infty$) the same as
\[h\times D({\bm \theta})\]
where we set
\[D({\bm \theta}):=\bigl \| \frac{\partial {\bf X}({\bm \theta})}{\partial \theta_1}
- \frac{\partial {\bf X}({\bm \theta})}{\partial \theta_2}\bigr\|_2.\]
Now, we may calculate
\begin{eqnarray*}
\frac{\partial X_1}{\partial\theta_1}-\frac{\partial X_1}{\partial \theta_2}&=&\frac{v'(\theta_1)(v^2(\theta_2)+v^2(\theta_3))+v'(\theta_2)v(\theta_1)v(\theta_2)}{W^{3/2}({\bm\theta})},\cr
\cr
\frac{\partial X_2}{\partial\theta_1}-\frac{\partial X_2}{\partial \theta_2}&=&-\frac{v'(\theta_2)(v^2(\theta_1)+v^2(\theta_3))+v'(\theta_1)v(\theta_1)v(\theta_2)}{W^{3/2}({\bm\theta})},\cr
\cr
\frac{\partial X_3}{\partial\theta_1}-\frac{\partial X_3}{\partial \theta_2}&=&-\frac{v(\theta_3)(v(\theta_1)v'(\theta_1)-v(\theta_2)v'(\theta_2))}{W^{3/2}({\bm\theta})}.
\end{eqnarray*}

$D^2$ is the sum of the squares of these three expressions. By elementary (but tedious) means it can be shown that
\[D^2({\bm \theta})=\Bigl(\frac{\pi}{2}\Bigr)^2
\frac{1+w(\theta_3)(1-w(\theta_1)-w(\theta_2))}
{(w(\theta_1)+w(\theta_2)+w(\theta_3))^2}.\]

\begin{lemma} We have the bounds. For $\theta_1+\theta_2+\theta_3=1,$ $\theta_j\ge0,$
\[1\le \frac{1+w(\theta_3)(1-w(\theta_1)-w(\theta_2))}
{(w(\theta_1)+w(\theta_2)+w(\theta_3))^2}\le  \frac{2\sqrt{3}+3}{3} \approx 2.15.\]
The lower bound is attained on the edge $\theta_3=0.$ 
\end{lemma}
\begin{proof}
For the lower bound just note that
\begin{align*}
\frac{1+w(\theta_3)(1-w(\theta_1)-w(\theta_2))}
{(w(\theta_1)+w(\theta_2)+w(\theta_3))^2}&\ge
\frac{1}
{(w(\theta_1)+w(\theta_2)+w(\theta_3))^2}\cr
&\ge 1
\end{align*}
as $w(\theta_1)+w(\theta_2)+w(\theta_3)\le1.$

For the upper bound, first note that the expression is of the form
\[F(w(\theta_1),w(\theta_2),w(\theta_3))\]
where $F(x_1,x_2,x_3)$ is symmetric in the first two variables, i.e., 
\[F(x_2,x_1,x_3)=F(x_1,x_2,x_3)).\]

To maximize it, on the boundary (where some $\theta_i=0$),  $w(\theta_1)+w(\theta_2)+w(\theta_3))=1$ and so
\begin{eqnarray*}F(w(\theta_1),w(\theta_2),w(\theta_3))&=&1+w(\theta_3)(1-w(\theta_1)-w(\theta_2))\cr
&\le& 1+w(\theta_3)\cr
&\le&2 <\frac{2\sqrt{3}+3}{3}.
\end{eqnarray*}

For a critical point, set $\theta_3=\-\theta_1-\theta_2$ so that
\[\frac{\partial}{\partial\theta_1} F(w(\theta_1),w(\theta_2),w(\theta_3))=
\frac{\partial F}{\partial x_1}w'(\theta_1)-\frac{\partial F}{\partial x_3}w'(\theta_3)\]
while
\begin{eqnarray*}
\frac{\partial}{\partial\theta_2} F(w(\theta_1),w(\theta_2),w(\theta_3))&=&
\frac{\partial}{\partial\theta_2} F(w(\theta_2),w(\theta_1),w(\theta_3))\cr
\cr
&=&\frac{\partial F}{\partial x_1}w'(\theta_2)-\frac{\partial F}{\partial x_3}w'(\theta_3).
\end{eqnarray*}
Hence critical points are given by
\begin{eqnarray*}
w'(\theta_1)&=&w'(\theta_3)\times\frac{\partial F}{\partial x_3}/\frac{\partial F}{\partial x_1},\cr
\cr
w'(\theta_1)&=&w'(\theta_3)\times \frac{\partial F}{\partial x_3}/\frac{\partial F}{\partial x_1}.
\end{eqnarray*}
In particular 
\[w'(\theta_1)=w'(\theta_2).\]
But this is possible  only if  either $\theta_2=1-\theta_1$ or $\theta_2=\theta_1.$ The first of these cases is a boundary point, and so we may assume that $\theta_2=\theta_1$ in which case
\[\theta_2=\theta_1=\frac{1-\theta_3}{2}.\]
Consequently,
\begin{eqnarray*}
w(\theta_2)=w(\theta_1)&=&w\Bigl(\frac{1-\theta_3}{2}\Bigr)\cr
\cr
&=&\frac{1-\cos\Bigl(\frac{\pi}{2}(1-\theta_3)\Bigr)}{2}\cr
\cr
&=&\frac{1-\sin\Bigl(\frac{\pi}{2}\theta_3\Bigr)}{2}\cr
\cr
&=&\frac{1-\sqrt{w(\theta_3)}}{2}.
\end{eqnarray*}
Setting $z:=\sqrt{w(\theta_3)}\in[0,1],$ the function to be maximized becomes
\begin{eqnarray*}
f(z)&:=&\frac{1+w(\theta_3)(1-w(\theta_1)-w(\theta_2))}{(w(\theta_1)+w(\theta_2)+w(\theta_3))^2}\cr
\cr
&=&\frac{1+z^2(1-2(1-z)/2)}{(2(1-z)/2+z^2)^2}\cr
\cr
&=&\frac{1+z^3}{(1-z+z^2)^2}\cr
\cr
&=&\frac{(1+z)(1-z+z^2)}{(1-z+z^2)^2}\cr
\cr
&=&\frac{1+z}{1-z+z^2},\quad z\in[0,1].
\end{eqnarray*}
Now,
\[f'(z)=\frac{2-2z-z^2}{(1-z+z^2)^2}\]
and so the critical point is for
\[z^2+2z-2=0\,\,\iff\,\,z=-1\pm\sqrt{3}.\]
The point $z=-1-\sqrt{3}$ is outside the interval $[0,1]$ and so the only critical point to be considered is $z=-1+\sqrt{3}.$ For this $z$ the value of $f$ is precisely
\[f(z)=\frac{2\sqrt{3}+3}{3}>2\]
and so this is the maximum value.
\end{proof}
\medskip
It follows that asymptotically the spacing of the spherical Waldron points are between 
\[\frac{1}{n}\times \frac{\pi}{2}\qquad {\rm and}\qquad\sqrt{\frac{2\sqrt{3}+3}{3}}\times \frac{1}{n}\times \frac{\pi}{2}.\]

\section{Polynomial Interpolation at the Waldron Points}

The quality of a set of interpolation points is perhaps best measured by the size of the Lebesgue constants, i.e.,
\[\Lambda_n:=\max_{\bf x\in S_d}\sum_{i=1}^N |\ell_i({\bf x})|.\]
Here we report some numerical results on the Lebesgue constants in dimensions two and three.

\subsection{Dimension Two}

We again consider the equilateral triangle centred at the origin with vertices
\[{\bf V}_1=(-\sqrt{3}/2,-1/2),\,\,{\bf V}_2=(+\sqrt{3}/2,-1/2),\,\,{\bf V}_3=(0,1).\]

For a comparison we also give the Lebesgue constants for the Simplex Points and also for the scheme based on concentric triangles introduced in \cite{Bos:83}.  These are also a (semi-)explicit set of points and to the best of our knowledge they are up to now the set of points for which there are semi-explicit formulas having the slowest order of growth of the Lebesgue constant. They are briefly described as follows (the details are provided in \cite{Bos:83} and also \cite{Bos:91}).

We place $3(n-3i)$ points on concentric triangles 
\[T_i\,\,\hbox{with}\,\,\hbox{vertices}\,\, R_i{\bf V_j},\,\, 1\le j\le 3\]
for $i=0,1,\cdots,s,$  radii 
\[1\ge R_0>R_1>\cdots>R_s>0,\] and $s$ defined by
\[s:=\Bigl\lfloor \frac{n-1}{3}\Bigr\rfloor.\]
In addition, if $n$ is a multiple of 3, then a single point is placed at the origin.  

The $3(n-3i)$ (of degree $m:=n-3i$) points on triangle $T_i$ are chosen to be the 3 vertices together with $m-1$ additional points spaced as Chebyshev-Lobatto points (i.e.,  $\cos(k\pi/m)$ for $[-1,1]$) on each edge, for a total of $3+3(m-1)=3m=3(n-3i.)$

The radii are chosen to maximize the associated Vandermonde determinant. Again, the details are explained in \cite{Bos:83}. It is this required maximization that makes the points only semi-explicit.

Here is a table of the radii, computed numerically, for degree up to 12.

\medskip

\begin{center}
\begin{tabular}{l|c|c}
      \textbf{Degree} & \textbf{No. Points} & \textbf{Radii}\\
      \hline
      1 & 3& 1\\
      2 & 6 & 1\\
      3 & 10 & 1,0\\
      4 & 15  & 1, \,$(1+3\sqrt{5})/22$\\
      5 & 21 & 1, \,0.5467133890977183\\
      6 & 28 & 1, \,0.6625914730317319,\,0\\
      7 & 36 & 1,\,0.7392097205159041,\,0.2099178922839476\\
      8 & 45 & 1,\,0.7926979593397175,\,0.3630731196442392\\
      9 & 55 & 1,\,0.8314018389721662,\,0.4713481792856927,\,0\\
      10 & 66 & 1,\,0.8603011832477779,\,0.5547886858166182,\\
      && 0.1489400918406532\\
      11 & 78 & 1,\,0.8824295392910452,\,0.6207291455415433,\\
      &&0.2691541556591404\\
      12 & 91 & 1,\, 0.8997282443826207,\,0.6734543809542708,\\
      && 0.3612491207621312,\,0
\end{tabular}
\end{center}

In Figure 3 we show the Concentric Triangle points (on the right) and the Waldron points, for degree $n=8.$

\begin{center}
 \includegraphics[width=15cm,height=10cm]{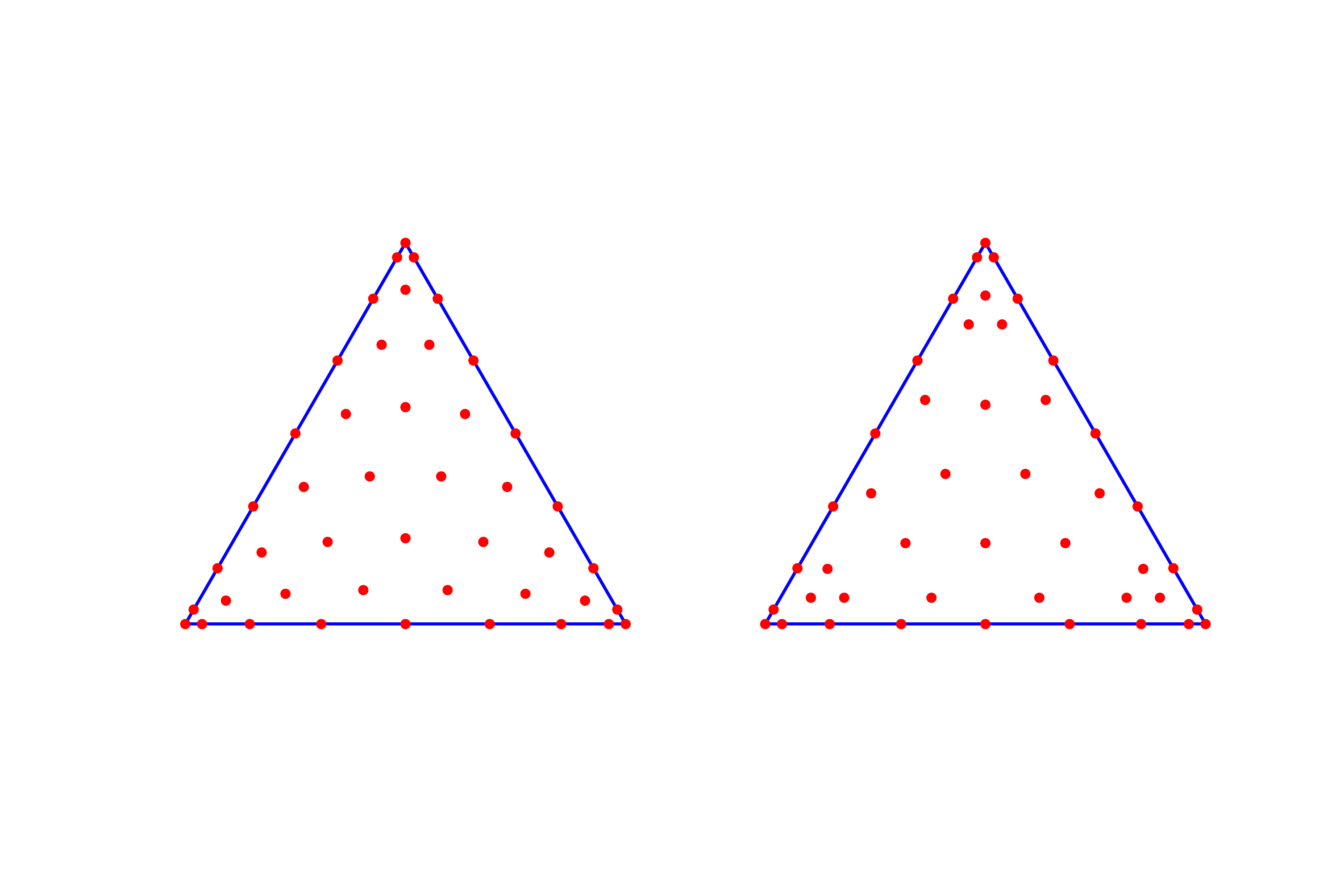}
\captionof{figure}{Degree 8 Waldron (left) and Concentric (right)}
\label{fig3}
\end{center}

\medskip

We now give the Lebesgue constants computed numerically by evaluating the Lebesgue function on a dense set of equally spaced points in the triangle.

\medskip 
\begin{center}
\begin{tabular}{|r|r|c|c|c|}
\hline
      $\bf n$  & $\bf N$& \textbf{Waldron Points}& \textbf{Concentric Points}&\textbf{Simplex Points} \\
      \hline
      1&3&1&1&1\\
      2&6&1.67&1.67&1.67\\
      3&10&2.11&2.11&2.27\\
      4&15&2.78&2.77&3.47\\
      5&21&3.36&4.11&5.45\\
      6&28&3.95&4.80&8.74\\
      7&36&4.63&6.01&14.34\\
      8&45&5.83&8.81&24.00\\
      9&55&7.18&10.75&40.87\\
      10&66&9.45&12.40&70.88\\
      11&78&12.37&18.28&124.52\\
      12&91&16.91&24.27&221.19\\
      13&105&23.34&$\cdot$&397.05\\
      14&120&33.04&$\cdot$&720.26\\
      15&136&47.38&$\cdot$&1315.77\\
      16&153&69.04&$\cdot$&2418.43\\
      \hline
\end{tabular}
\end{center}

\subsection{Dimension Three}

For dimension 3 the Waldron points per se do not have low growth of the Lebesgue constant. However, if with a slight variant we do obtain somewhat better growth.  This variant is to force the points on a face of the simplex to be the two dimensional Waldron points. Specifically, in Definition \ref{waldronpts}, we set $\omega_j:=0$ if $\alpha_j=0$  and ensure that the sum of the other three coordinates is one.  Here we compare the Lebesgue constants for the modified Waldron points and the Simplex points.

\medskip 
\begin{center}
\begin{tabular}{|r|r|c|c|}
\hline
      $\bf n$  & $\bf N$& \textbf{Waldron Points}&\textbf{Simplex Points} \\
      \hline
      1&4&1&1\\
      2&10&2.00&2.00\\
      3&20&2.99&3.02\\
      4&35&4.25&4.89\\
      5&56&5.49&8.08\\
      6&84&7.68&13.65\\
      7&120&10.15&23.37\\
      8&165&14.57&40.45\\
      9&220&21.06&71.00\\
      10&286&33.00&126.13\\
      11&364&56.00&225.42\\
      12&455&90.63&406.01\\
      \hline
\end{tabular}
\end{center}


\end{document}